\documentclass{rspublic}[]

\usepackage{graphicx}

\renewcommand\({\left(}
\renewcommand\){\right)}

\newcommand{\e}{{\rm e}}
\newcommand{\half}{{\mbox{$\frac{1}{2}$}}}

\renewcommand{\[}{\begin{equation}}
\renewcommand{\]}{\end{equation}}

\newcommand{\OU}{Ornstein--Uhlenbeck}

\newlength{\figurewidth}
\let\epsilon=\varepsilon

\newcommand{\Var}{{\mbox{Var}}}

\begin{document}

   \title{Stochastic Stokes' drift with inertia}

   \author[K. M. Jansons]{Kalvis M. Jansons%
   \footnote{Email: \texttt{stokesdrift@kalvis.com}}}
   \affiliation{Department of Mathematics,
   University College London,
   Gower Street,
   London WC1E 6BT, UK}
   \label{firstpage}

   \maketitle

   \begin{abstract}{Brownian motion, Stochastic Stokes' drift,
       particle sorting}{We consider both the effect of particle inertia on
       stochastic Stokes' drift, and also a related process
       which could be considered as a crude model of stochastic
       Stokes' drift driven by an eddy diffusivity.  In the latter,
       the stochastic forcing is a stable \OU{} process rather than
       Brownian motion.  We show that the eddy Stokes' drift velocity
       has a peak at a non-zero value of the correlation time-scale
       for particles that have the same (limiting) diffusivity.  For
       both of the models considered, this study shows that not only
       can stochastic Stokes' drift be used to sort particles with
       different diffusivities, but also it can be used to sort
       particles of the same diffusivities but with different particle
       masses or correlation time-scales.  This effect may be
       important in particle sorting applications.}
   \end{abstract}

   \InputIfFileExists{KalvjaOpts.tex}

   \section{Introduction}
   \label{intro}
   
   This study is a continuation of the work of Jansons \& Lythe (1998)
   on stochastic Stokes' drift, which is the modification of the
   classical Stokes' drift of a particle in a travelling wave due to
   the effect of Brownian motion (or some other random forcing with
   very short timescale correlations). In Jansons \& Lythe (1998), the
   particle inertia was neglected.  In this study, we first consider,
   in \S\ref{inertia}, the effect of particle inertia on the
   stochastic Stokes' drift velocity.  For simplicity, we still
   neglect the inertia in the fluid, so we can use the Stokes flow
   approximation to determine the drag on the particle, and so that
   this drag is not history dependent. This is a consistent
   approximation for a particle in a gas, where it is reasonable to
   assume that the particle density is much greater than that of the
   ambient fluid.  For the motion of a particle in a liquid of similar
   density to that of the particle, the situation is much more
   complex, and the work here is not applicable, though may agree
   qualitatively in some cases.
   
   Many authors have continued the investigation of various aspects of
   stochastic Stokes' drift.  The first exact solution to Brownian
   motion forced Stokes' drift is due to Van Den Broeck (1999), who
   found an exact solution for an arbitrary wave in an integral form
   and simple exact solutions for a square wave and an impulse wave.
   Applications of stochastic Stokes' drift to ocean circulation have
   recently been considered by Restrepo (2006), and the techniques
   considered in this study appear to be relevant.  Li Yu-Xiao et. al.
   (2001) studied the effect of asymmetric potentials on the
   stochastic Stokes' drift. Bena et. al. (2005) considered, both
   analytically and numerically, Stokes drift forced by a dichotomous
   Markov process.
   
   Also the theory of stochastic Stokes' drift applies to Brownian
   motors and tilting ratchets (see Reimann (2002)).

   Jansons \& Lythe (1998) showed that particles of different diffusivities
   could be separated by the superposition of several waves with
   different directions, wavenumbers and frequencies.  In such a
   system, it is possible to arrange for particles to drift in
   different directions depending on their diffusivities.  We show
   here that it is also possible to separate particles of the same
   diffusivity but with different masses.  It is not yet clear if this
   effect will have practical applications, but the early signs are
   promising.
   
   In \S\ref{eddy}, we consider a related process of independent
   interest, and which could be considered a crude model of a particle
   moving in a travelling wave and small eddies, i.e. the random
   motion has short time-scale correlations.  This is mathematically
   stochastic Stokes' drift forced not by Brownian motion but rather
   by a stable \OU{} process.  We show, for example, that it is
   possible for the stochastic Stokes' drift velocity to be increased
   by a careful choice of the correlation time-scale, even for
   particles with the same diffusivity in the long-time limit.  Such
   an effect might be important in oceanographic applications of
   particle dispersion for example.

   In \S\ref{results}, we compare the leading-order asymptotic
   approximations of the Stokes' drift for the two types of models
   considered here with Monte Carlo simulations, which agree well even
   for reasonably large values of the `small' asymptotic parameter.

   Finally, in \S\ref{morewaves}, we consider applications of these
   results to particle separation schemes using multiple waves in
   dimensions greater than one.

   \section{Stochastic Stokes' drift with inertia}
   \label{inertia}
   
   As shown by Jansons \& Lythe (1998), stochastic Stokes' drift in
   higher dimensions with many waves is, to leading order, a
   superposition of the contributions from each wave separately.  This
   is because, in the leading-order calculation, cross-terms from
   different waves average to zero if they have either a different
   spatial or different temporal frequency.  This is also true for the
   generalizations considered in this study.  So, for simplicity, we
   consider the essentially one-dimensional system with a single wave,
   and comment on the higher-dimensional versions in \S\ref{morewaves}.

   The standard form of the Langevin equation for the velocity $U$ of
   an isotropic Brownian particle in a gas at rest is given by
   \[\label{langevin}
   m \frac{dU}{dt} = - b^{-1} U + \eta,
   \]
   where $m$ and $b$ are respectively the particle mass and particle
   mobility, both assumed constant, and $\eta$ is a time-dependent
   noise term.  
   
   Note that assuming that $b$ is a constant follows from the
   assumption that we can neglect inertia effects in the motion of the
   fluid around the particle.  This is a rational approximation in the
   case of a particle in a fluid, provided that the particle density
   is much greater than that of the fluid, which is normally so for a
   gas, but is rarely true for a liquid.  In the case of a liquid, the
   particle drag term depends on the history of the particle's
   motion (see Hinch (1975)), and we need to include both fluid and
   particle inertia terms in the description of the motion, which we
   shall not do here.  However, the results for a particle in a liquid
   are likely to be qualitatively similar to those of this study.

   Turning \eqref{langevin} into an It\^o stochastic differential
   equation we find
   \[\label{ItoLangevin}
   dU_t = -\lambda U_t dt + \lambda \sigma dB_t,
   \]
   where $B$ is a standard Brownian motion (i.e.\ a Wiener process),
   $\lambda = (b m)^{-1}$, $\sigma = (2 b K {\cal T})^{1/2}$, $K$ is
   Boltzmann's constant, and ${\cal T}$ is the (absolute) temperature.
   The long-time particle diffusivity is $\half \sigma^2$.

   We now introduce an additional velocity term $\epsilon f$ due
   to some wave motion imposed on the fluid, giving
   \[
   dU_t = -\lambda (U_t - \epsilon f(X_t, t)) dt + \lambda \sigma dB_t,
   \]
   where $X_t$ is the particle position and $\epsilon$ is a
   dimensionless constant, which will be used as a small parameter in
   the asymptotic analysis below.  The stochastic Stokes' drift
   velocity is defined as
   \[
   V \equiv \lim_{t\to\infty} t^{-1}\(X_t - X_0\).
   \]
   Using the Ergodic theorem we are able to replace this limit with an
   expectation of $U_t$, provided the starting phase of the wave is
   uniformly distributed.

   For a wave with characteristic wavenumber $k$ and angular frequency
   $\omega$, the natural time-scales in this system are
   \[\label{timescales}
   \lambda^{-1}, \qquad \omega^{-1}, \qquad \sigma^{-2} k^{-2},
   \]
   which are respectively the relaxation time-scale of the particle
   velocity, the period of the wave and the time for the particle to
   diffuse over a wavelength.  These time-scales give two independent
   non-dimensional constants, in addition to $\epsilon f_0 k /\omega$,
   where $f_0$ is a typical value of $f$.  We assume that all of the
   time-scales in \eqref{timescales} are comparable, and denote the
   common time-scale as $T$.
 
   We can expand both $U$ and $X$ as a formal asymptotic series in
   $\epsilon$, in the limit $\epsilon\to0$, namely
   \[
   U_t = U^{(0)}_t + \epsilon U^{(1)}_t + \epsilon^2 U^{(2)}_t + \cdots
   \]
   and 
   \[
   X_t = X^{(0)}_t + \epsilon X^{(1)}_t + \epsilon^2 X^{(2)}_t + \cdots
   \]
   These local approximations are valid for an $o(\epsilon^{-1} T)$
   time range, which is sufficient for the determination of stochastic
   Stokes' drift as the exponential tails beyond this time range do
   not contribute to the leading-order asymptotic results.  At times
   of order $\epsilon^{-1} T$ some contributions switch order, but we
   do not need to address this point as we are not going to determine
   the first correction to the stochastic Stokes' drift in this study.

   We now consider the formal $\epsilon$ expansion an order at a time.
   
   \subsection{Order $\epsilon^0$}
   \[
   \lambda^{-1} dU^{(0)}_t = -U^{(0)}_t dt + \sigma dB_t,
   \]
   which is a stable \OU{} process and can be solved exactly in the
   form
   \[\label{U0}
   U^{(0)}_t = \e^{-\lambda t} W\(\half\lambda\sigma^2 \e^{2\lambda t}\),
   \]
   where $W$ is a standard Brownian motion, with $W(0)=0$, which means
   that we have started $U^{(0)}$ with its stationary law.  Note that
   $U^{(0)}_t$ has a zero mean.

   Thus the autocorrelation of $U^{(0)}$ is given by
   \[
   E[U^{(0)}_s U^{(0)}_t] = \half \lambda \sigma^2 \exp(-\lambda |s-t|).
   \]
   
   \subsection{Order $\epsilon^1$}
   \[
   \lambda^{-1} \frac{dU^{(1)}_t}{dt} = -U^{(1)}_t + f(X^{(0)}_t, t).
   \]
   Thus 
   \[
   U^{(1)}_t = \lambda \int_0^\infty 
   \exp(-\lambda \tau) f(X^{(0)}_{t-\tau}, t-\tau)
   d\tau.
   \]

   \subsection{Order $\epsilon^2$}
   It is at this order we find the first contribution to the
   stochastic Stokes' drift.
   \[
   \lambda^{-1} \frac{dU^{(2)}_t}{dt} =  -U^{(2)}_t
   + f'(X^{(0)}_t, t) X^{(1)}_t,
   \]
   where the dash denotes differentiation with respect to the first argument.
    Thus
   \[
   \frac{U^{(2)}_t}\lambda 
   = \int_0^\infty
   \exp(-\lambda \tau) f'(X^{(0)}_{t-\tau}, t-\tau) X^{(1)}_{t-\tau}
   d\tau.
   \]
   This gives
   \[
   \frac{U^{(2)}_t}\lambda 
   = \int_0^\infty \!\int_\tau^\infty
   \exp(-\lambda \tau) f'(X^{(0)}_{t-\tau}, t-\tau)
   U^{(1)}_{t-\alpha} 
   d\alpha d\tau,
   \]
   and so
   \[\label{U2}
   \frac{U^{(2)}_t}{\lambda^2} 
   = \int_0^\infty \!\int_\tau^\infty \!\int_0^\infty
   \exp(-\lambda (\tau + \beta)) f'(X^{(0)}_{t-\tau}, t-\tau)
   f(X^{(0)}_{t-\alpha-\beta}, t-\alpha-\beta)
   d\beta d\alpha d\tau.
   \]
   This is about as far as it is useful to go in general, so we now
   focus on the particular case
   \[
   f(x,t) = u \cos(k x - \omega t + \phi),
   \]
   where $u$, $k$, $\omega$ and $\phi$ are constants. The phase $\phi$
   of the wave will not alter the long-time limits of interest, so it
   is not important. However, for convenience, and to avoid starting
   transients, we choose $\phi$ to have a uniform distribution on
   $[0,2\pi)$.  Note that in this special case, $E[U^{(0)}_t] =
   E[U^{(1)}_t] = 0$.
     
   So the leading-order Stokes' drift is from $U^{(2)}$. From
   \eqref{U2}, we find
   \[\label{U2S}
   \frac{U^{(2)}_t}{\lambda^2 u^2 k} 
   = - \int_0^\infty \!\int_\tau^\infty \!\int_0^\infty
   \exp(-\lambda (\tau + \beta)) \sin(\phi_s) \cos(\phi_c)
   d\beta d\alpha d\tau,
   \]
   where $\phi_s = k X^{(0)}_{t-\tau} - \omega (t-\tau) + \phi$ and
   $\phi_c = k X^{(0)}_{t-\alpha-\beta} - \omega (t-\alpha-\beta) + \phi$.
   On averaging we find
   \[\label{meanU2}
   \frac{E[U^{(2)}_t]}{\half \lambda^2 u^2 k} 
   = - \int_0^\infty \!\int_\tau^\infty \!\int_0^\infty
   \exp(-\lambda (\tau + \beta)) E[\sin(\phi_s - \phi_c)]
   d\beta d\alpha d\tau.
   \]
   This is exact, rather than a long-time limit, since we chose $\phi$ to
   have a uniform distribution, so the time $t=0$ is not special.
   
   Note that
   \[
   \phi_s - \phi_c = k (X^{(0)}_{t-\tau} - X^{(0)}_{t-\alpha-\beta})
   + \omega(\tau-\alpha-\beta),
   \]
   which is Gaussian, and so $E[\sin(\phi_s - \phi_c)]$ can be
   evaluated exactly.  From \eqref{U0}, we see that
   \[
   C(t-s) \equiv E[(X^0_t-X^0_s)^2] = 
   \sigma^2 \(|t-s| + \frac{(\e^{-\lambda |t-s|} - 1)}\lambda\).
   \]
   Thus
   \[
   E[\sin(\phi_s - \phi_c)]
   = \exp(-\half k^2 C(\tau-\alpha-\beta)) \sin(\omega(\tau-\alpha-\beta)).
   \]
   Putting everything together, we find from \eqref{meanU2} the
   stochastic Stokes drift velocity, $V$, to leading-order in $\epsilon$:
   \[
   \frac{V}{\half \epsilon^2 \lambda^2 u^2 k} 
   = \int_0^\infty \!\int_\tau^\infty \!\int_0^\infty
   \exp(-\lambda (\tau + \beta) - \half k^2 C(\alpha + \beta - \tau)) 
   \sin(\omega(\alpha + \beta - \tau))
   d\beta d\alpha d\tau.
   \]
   We may simplify this a little, to give
   \[\label{inertiadrift}
   V
   = {\half \epsilon^2 \lambda u^2 k} 
   \int_0^\infty \!\int_0^\infty
   \exp(-\lambda \beta - \half k^2 C(\alpha + \beta)) 
   \sin(\omega(\alpha + \beta))
   d\beta d\alpha.
   \]
   
   We should expect \eqref{inertiadrift} to reduce to the result of
   Jansons \& Lythe (1998) in the limit $\lambda\to\infty$, i.e.\ in
   the zero mass limit.  Taking this limit reduces \eqref{inertiadrift}
   to 
   \[\label{ilimit1}
   V
   = {\half \epsilon^2 u^2 k}
   \int_0^\infty
   \exp(- \half k^2 C(\alpha))
   \sin(\omega\alpha)
   d\alpha,
   \]
   and
   \[\label{ilimit2}
   C(\alpha) = \sigma^2 |\alpha|.
   \]
   Equations \eqref{ilimit1} and \eqref{ilimit2} agree exactly with
   the penultimate line in the calculation of Jansons \& Lythe (1998),
   but in a slightly different notation.

   To leading order
   \[
   \lim_{t\to\infty} \frac{\Var[X_t]}t = \sigma^2.
   \]
   We could easily determine the first correction to the variance of
   position, but we do not present this result here as it was felt
   that it was a little too messy to be of interest.

   \section{Eddy Stokes' drift}
   \label{eddy}
   
   In this section, we consider a closely-related process to the
   stochastic Stokes' drift with inertia process of \S\ref{inertia}.
   This process is stochastic Stokes' drift forced by a stable \OU{}
   process rather than Brownian motion.  The governing equation (using
   notation as close to \S\ref{inertia} as possible) is
   \[
   \frac{dX_t}{dt} = U_t + \epsilon f(X_t, t),
   \]
   where $U_t$ satisfies \eqref{ItoLangevin}, but, of course, does not
   have the same physical interpretation.  This process could be
   considered as a crude (or qualitative) description of stochastic
   Stokes' drift with the randomness coming from eddy diffusivity.
   
   In the case $f(x,t) = u \cos(k x - \omega t + \phi)$, the
   stochastic Stokes' drift, to leading order, is given by
   \[
   V = \half \epsilon^2 u^2 k \int_0^\infty \sin(\omega t)
   \exp\(-\half k^2 C(t)\) dt,
   \]
   which follows from an argument that is a simpler version of that in
   \S\ref{inertia}, and is the same as that of Jansons \& Lythe (1998)
   but with $\sigma^2 t$ replaced by $C(t)$, so will not be given
   here.  In this case, it is also trival to see that the
   $\lambda\to\infty$ limit of this result is exactly the same as that
   for Brownian forcing in Jansons \& Lythe (1998).
   
   We now determine the leading-order correction to the long-time
   limit of the variance of position.  This will appear at order
   $\epsilon^2$. We take $X^{(n)}_0 = 0$ for all $n$.  Note that, to the
   required order, $t^{-1}E[X_t]\to 0$, as $t\to\infty$, so we need
   only to determine
   \[\label{varX}
   E[X_t^2] = E\left[(X_t^{(0)})^2\right] 
   + 2\epsilon E\left[X_t^{(0)}X_t^{(1)}\right]
   + \epsilon^2\(E\left[(X_t^{(1)})^2\right] 
   + 2 E\left[X_t^{(0)}X_t^{(2)}\right]\)
   + \cdots,
   \]
   which we consider term by term.
   \[
    t^{-1} E\left[(X_t^{(0)})^2\right] 
    = t^{-1} C(t)
    \to \sigma^2,
    \]
    as $t\to\infty$;
    \[
    t^{-1} E\left[X_t^{(0)}X_t^{(1)}\right] \to 0,
    \]
    as $t\to\infty$;
    \[
    \begin{aligned}
      t^{-1} E\left[(X_t^{(1)})^2\right] 
      &= t^{-1} u^2 \int_0^t \int_0^t 
      E\left[\cos\(k X_\alpha^{(0)} - \omega \alpha\)
        \cos\(k X_\beta^{(0)} - \omega \beta\)\right]
      d\alpha d\beta\\
      &\to u^2 \int_0^\infty \cos(\omega t) 
      \exp\(-\half k^2 C(t)\) dt,
    \end{aligned}
    \]
    as $t\to\infty$;
    \[
    t^{-1} E\left[X_t^{(0)}X_t^{(2)}\right] \to 0,
    \]
    as $t\to\infty$.
    
    Thus
    \[\label{EddyVar}
    \lim_{t\to\infty}\frac{\Var[X_t]}t
    =
    \sigma^2 
    + \epsilon^2 u^2 
    \int_0^\infty \cos(\omega t) \exp\(-\half k^2 C(t)\) dt
    + \cdots.
    \]

    We can also consider the limit $\lambda\to\infty$ of this result,
    giving the corresponding result for Brownian forcing, though this
    was not determined by Jansons \& Lythe (1998) in their study of
    this process.  In this limit, we find
    \[
    \lim_{t\to\infty}\frac{\Var[X_t]}t
    =
    \sigma^2
    + \epsilon^2 \frac{2 u^2 k^2 \sigma^2}{k^4 \sigma^4 + 4 \omega^2}
    + \cdots.
    \]

   \section{Numerical results}
   \label{results}
   
   To determine how small $\epsilon$ needs to be for the asymptotic
   results to be a good approximation to the exact results, we perform
   some Monte Carlo simulations. In the simulations shown here, we use
   the stochastic Euler method with $dt = 0.001$ and all parameters
   other than $\epsilon$ set equal to $1$.  These simulations are run
   to $t=10^8$ to control the variance sufficiently for comparison of
   the time-averaged stochastic Stokes' drift.

   We present numerical confirmation of the results of both
   \S\ref{inertia} and \S\ref{eddy}.
   
\begin{figure}
     \centering
     \includegraphics[angle=-90,width=\figurewidth]{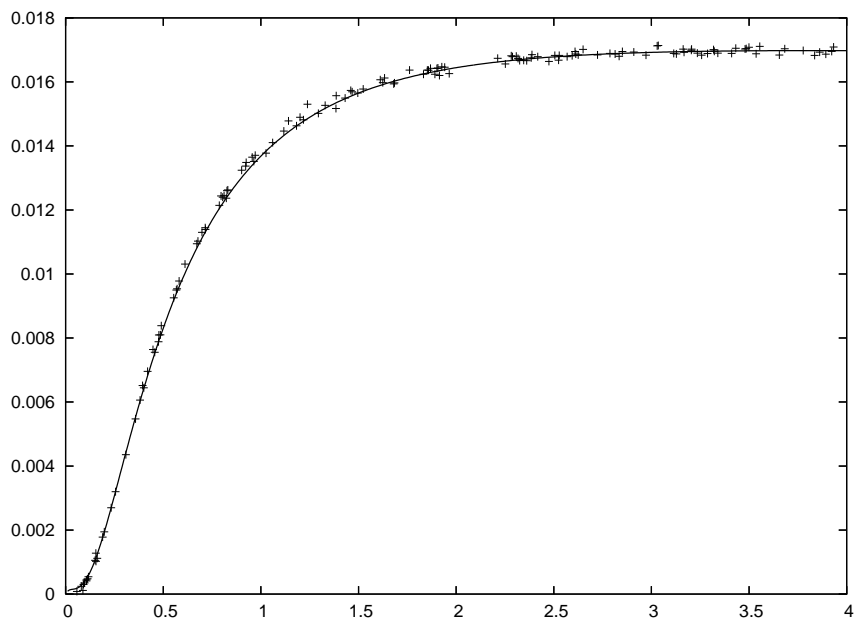}
     \caption{Stochastic Stokes' drift velocity against $\lambda$ for
     $\epsilon=0.2$, and all other parameters set equal to $1$. The
     curve is the leading-order asymptotic result and the points are
     Monte Carlo results.}
     \label{inertia1}
   \end{figure}

   In the case of stochastic Stokes' drift with inertia, and for
   $\epsilon=0.2$, the agreement is excellent and is shown in Figure
   \ref{inertia1}, though the (predicted) variance leads to a little
   scatter around the asymptotic result.  This could have been reduced
   by running the Monte Carlo simulations up to a a larger time $t$,
   but as it was each point took $5$ hours on a $3\mbox{GHz}$ Intel
   CPU.
   
\begin{figure}
     \centering
     \includegraphics[angle=-90,width=\figurewidth]{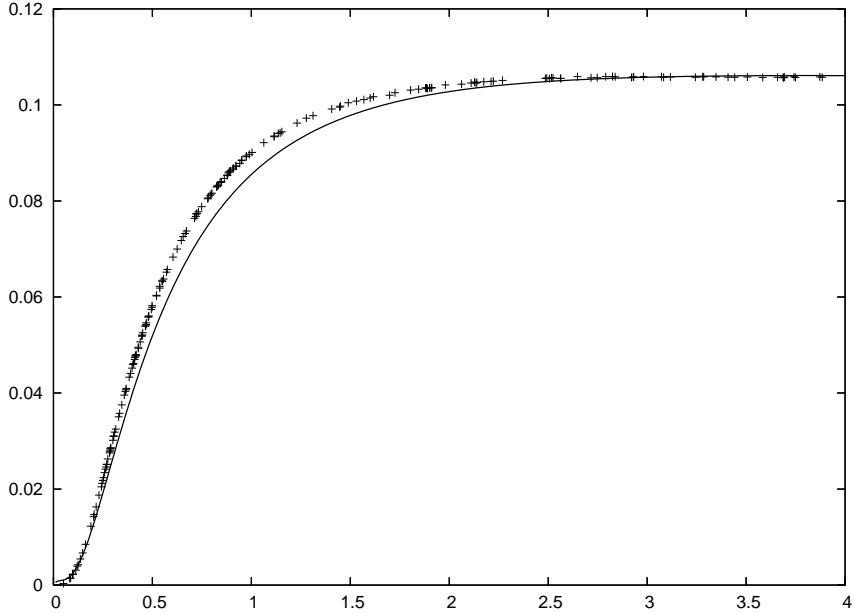}
     \caption{Stochastic Stokes' drift velocity against $\lambda$ for
     $\epsilon=0.5$, and all other parameters set equal to $1$. The
     curve is the leading-order asymptotic result and the points are
     Monte Carlo results.}
     \label{inertia2}
   \end{figure}

   When $\epsilon=0.5$ for this model, the scatter of the Monte Carlo
   results is reduced (see Figure \ref{inertia2}) and a small
   departure from the asymptotic approximation can be seen, but it is
   still small enough for the asymptotic result to be useful in
   applications.

\begin{figure}
     \centering
     \includegraphics[angle=-90,width=\figurewidth]{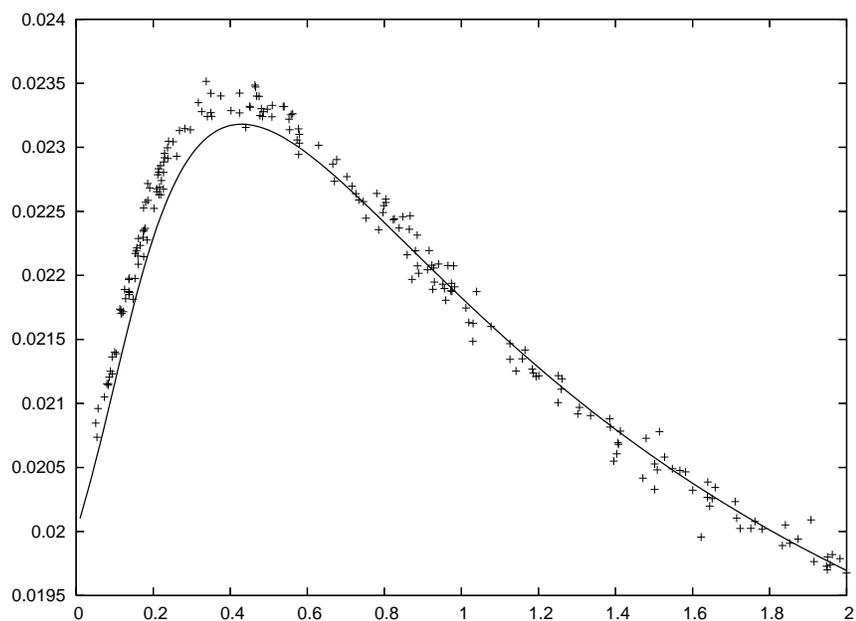}
     \caption{Eddy Stokes' drift velocity against $\lambda$ for
     $\epsilon=0.2$, and all other parameters set equal to $1$. The
     curve is the leading-order asymptotic result and the points are
     Monte Carlo results.}
     \label{eddy1}
   \end{figure}

\begin{figure}
     \centering
     \includegraphics[angle=-90,width=\figurewidth]{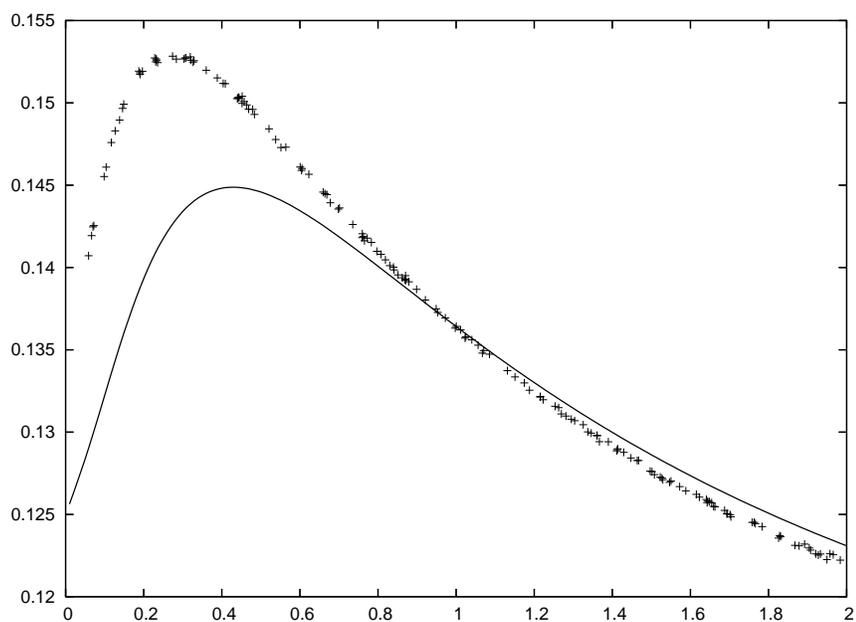}
     \caption{Eddy Stokes' drift velocity against $\lambda$ for
     $\epsilon=0.5$, and all other parameters set equal to $1$. The
     curve is the leading-order asymptotic result and the points are
     Monte Carlo results.}
     \label{eddy2}
   \end{figure}

   The agreement between Monte Carlo results and asymptotic results is
   much the same for the case of eddy Stokes' drift (\S\ref{eddy}).
   However, for this model there is a peak in the Stokes' drift, which
   was not seen in the model of \S\ref{inertia}, and is shown in
   Figures \ref{eddy1} and \ref{eddy2} for $\epsilon=0.2$ and
   $\epsilon=0.5$ respectively.  Thus the eddy Stokes' drift at finite
   $\lambda$ can be larger than that of the $\lambda\to\infty$ limit
   (i.e. larger than for the Brownian Stokes' drift result).  This
   does not appear to occur in the model of \S\ref{inertia}.

\begin{figure}
     \centering
     \includegraphics[angle=-90,width=\figurewidth]{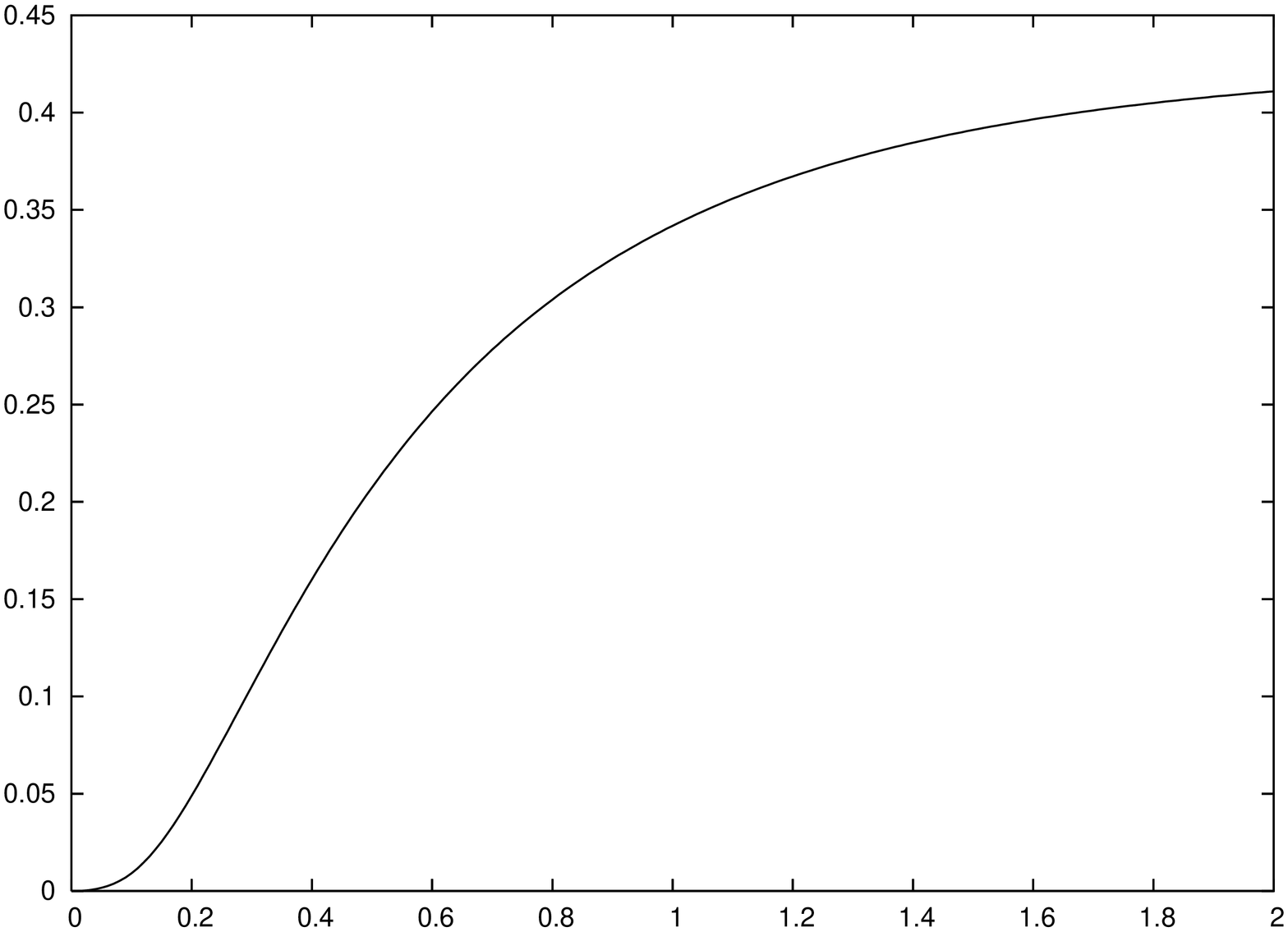}
     \caption{The first correction to the variance coefficient for the
     particle position for eddy Stokes' drift against $\lambda$.}
     \label{eddy3}
   \end{figure}

   Finally, we plot the first correction to the large-time limit of the
   particle variance for eddy Stokes' drift in Figure \ref{eddy3}.

   \section{Multiple waves and higher dimensional results}
   \label{morewaves}
   
   In the work of Jansons \& Lythe (1998), they considered stochastic
   Stokes' drift in several dimensions. In this case, it was possible to
   sort particles by a combination of waves in different directions
   and different spatial or temporal frequencies.  At leading order in
   the Stokes' drift, such combinations of waves act independently on
   the particle as the cross-terms average to zero in the long-time
   limit. Thus the modifications to the classical Stokes' drift due to
   diffusion is different for each of the wave components and so the
   linear combinations of these drifts will in general result in
   particles of different diffusivities diffusing in different
   directions.  This leads to an interesting continuous particle
   sorting method.
   
   In the current study too, the higher dimensional versions work in
   exactly the same way, as for the same reason all cross-terms due to
   different waves average to zero, at leading order, provided that
   they have either different spatial or temporal frequencies.  Thus
   again the resulting leading-order stochastic Stokes' drift is just
   a linear combination of the stochastic Stokes' drifts of each
   component.

   This implies that, in the case of stochastic Stokes' drift with
   inertia (see \S\ref{inertia}), it is possible to sort particles that
   have the same diffusivity but different masses, by arranging them
   to have different directions for the resulting Stokes' drift.  Whether
   this has practical applications remains to be seen, but it certainly
   looks hopeful.

   Equally, in the case of eddy Stokes' drift (see \S\ref{eddy}),
   particles with the same diffusivity in the long-time limit but with
   different correlation times can be sorted by arranging the
   resulting Stokes' drift to be in different directions.

   \section{Conclusions}

   In this study, we have extended the work of Jansons \& Lythe (1998) in
   two directions. Firstly in \S\ref{inertia}, we consider stochastic
   Stokes' drift including the effect of particle inertia, and
   secondly, we consider stochastic Stokes' drift forced by a stable
   \OU{} process rather than Brownian motion.  Both of these results
   are to leading-order in the strength $\epsilon$ of the wave motion,
   and the stochastic Stokes' drift appears at order $\epsilon^2$.
   These results extend trivially to higher dimensions and many waves,
   as the cross-terms for waves of different spatial or temporal
   frequencies average to zero.  Thus the `fanout' of particles of
   different diffusivities observed by Jansons \& Lythe (1998)
   extends to particles of the same (long-time limit) of diffusivity
   but with different masses in the case of the model in
   \S\ref{inertia} and different correlation times in the case of the
   model of \S\ref{eddy}.  In the case of eddy Stokes' drift (see
   \S\ref{eddy}), we also compute the long-time limit of the variance
   of particle positions.  The short correlation time limit of this
   result applies directly to the Brownian motion forcing case
   studied by Jansons \& Lythe (1998).
   
   The agreement of these results with Monte Carlo simulations is
   impressive (see \S\ref{results}), and shows that, for practical
   applications, the leading-order asymptotic approximations are
   sufficient, even when the `small' asymptotic parameter $\epsilon =
   0.5$.  Note that this is partly due to the asymptotic expansion
   being naturally in $\epsilon^2$ rather than $\epsilon$, as
   replacing $\epsilon$ with $-\epsilon$ is equivalent to a phase
   shift, and does not effect the stochastic Stokes' drift velocity.
   
   It appears that Van Den Broeck's (1999) exact method could be extended
   to consider the system of this study, but it is less clear if there
   are simple exact expressions even for a square wave.
   
   One extension of these results, which would be of interest, is to
   consider the behaviour of flexible particles that are large enough
   to feel different parts of the wave.  The simplest such example
   would be a flexible dumbbell, consisting of two point particles
   connected by a linear spring.

   \begin{verbatim}
     $Id: StokesDrift.tex 2343 2006-09-06 12:44:57Z kalvis $
   \end{verbatim}
   
\end{document}